  \documentclass[letterpaper, 10pt, twocolumn]{ieeeconf}

\usepackage{epsfig} %

\usepackage{amsmath,graphicx,amsfonts,amssymb,epsfig,subfig,mathrsfs,mathtools}
\usepackage{ntheorem}
\usepackage{cuted}
\usepackage{arydshln}
\usepackage{blkarray}

\usepackage{color}
\usepackage{multirow}
\usepackage{multicol}
\usepackage{lipsum}
\usepackage{rotating}
\usepackage{graphicx}%
\usepackage{psfrag}
\usepackage{algorithm}

\usepackage[utf8]{inputenc}
\usepackage{xspace}

\usepackage{algorithmicx}
\usepackage{algpseudocode}
\usepackage{cite}
\usepackage{tikz}
\usepackage{cancel}
\usepackage{textcomp}
 \definecolor{lin}{RGB}{240,0,0}
 \definecolor{paleblue}{RGB}{0,9,255}
\usepackage[colorlinks=true,linkcolor=lin,citecolor=paleblue,pdfa]{hyperref}
\usepackage{comment}
\includecomment{comment}

\newcommand{\maps}[3]{#1: #2 \mapsto #3}
\newcommand{\map}[3]{#1: #2 \rightarrow #3}
\newcommand{\setdef}[2]{\{#1 \; | \; #2\}}
\newcommand{\st}{\ensuremath{\operatorname{s.t.}}}
\newcommand{\real}{\ensuremath{\mathbb{R}}}
\newcommand{\prob}{\ensuremath{\mathbb{P}}}

\newcommand{\realnonnegative}{\ensuremath{\mathbb{R}}_{\ge 0}}

\newcommand{\supscr}[2]{#1^{\textup{#2}}}
\newcommand{\vect}[1]{\boldsymbol{#1}}

\newcommand{\Norm}[1]{\|#1\|}
\newcommand{\trans}[1]{{#1}^\top}

\newcommand{\Prob}{\ensuremath{\operatorname{Prob}}}

\newcommand{\untilinterval}[2]{\{{#1},\dots, {#2}\}}

\newcommand{\Batch}{\ensuremath{\operatorname{B}}}
\newcommand{\Detect}{\ensuremath{\operatorname{D}}}
\newcommand{\Alarm}{\ensuremath{\operatorname{Alarm}}}
\newcommand{\attack}{\ensuremath{\operatorname{attack}}}
\newcommand{\noattack}{\ensuremath{\operatorname{no \; attack}}}

\newcommand{\reviseone}[1]{{#1}}

\DeclarePairedDelimiterX\set[1]\lbrace\rbrace{\setaux#1}
 \def\setaux#1|{#1\;\delimsize\vert\;}

\makeatletter
\newtheoremstyle{breaknote}%
  {\item{\theorem@headerfont
          ##1\ ##2\theorem@separator}\hskip\labelsep\relax}%
  {\item{\theorem@headerfont
          ##1\ ##2\ (##3)\theorem@separator}\hskip\labelsep\relax}
\makeatother
\theoremstyle{breaknote}
\newtheorem{assumption}{Assumption}[section]
\newtheorem{theorem}{Theorem}[section]

\newtheorem{definition}{Definition}[section]
\newtheorem{lemma}{Lemma}[section]

\newtheorem{remark}{Remark}[section]

\allowdisplaybreaks

\title{High-Confidence Attack Detection via Wasserstein-Metric
  Computations}

\author{Dan Li$^{1}$ and Sonia Mart{\'\i}nez$^{1}$
  \thanks{This research was developed with funding from
    ONR N00014-19-1-2471, and AFOSR FA9550-19-1-0235.}
  \thanks{$^{1}$D. Li and S. Mart{\'\i}nez are with the Department of
    Mechanical and Aerospace Engineering, University of California San
    Diego, La Jolla, CA 92092, USA {\tt\small lidan@ucsd.edu;
      soniamd@ucsd.edu}}%
}

\begin{document}

\maketitle

\begin{abstract}
  This paper considers a sensor attack and fault detection problem for
  linear cyber-physical systems, which are subject to
  \reviseone{system} noise that can \reviseone{obey an} unknown
  light-tailed distribution. We propose a new threshold-based
  detection mechanism that employs the Wasserstein metric, and which
  guarantees system performance with high confidence
  \reviseone{employing a finite number of measurements}. The proposed
  detector may generate false alarms with a rate $\Delta$ in normal
  operation, where $\Delta$ can be tuned to be arbitrarily small by
  means of a \textit{benchmark distribution} which is part of our
  mechanism.  Thus, the proposed detector is sensitive to sensor
  attacks and faults which have a statistical behavior that is
  different from that of the system noise. We quantify the impact of
  \textit{stealthy} attacks---which aim to perturb the system
  operation while producing false alarms that are consistent with the
  natural system noise---via a \textit{probabilistic} reachable
  set. To enable tractable implementation of our methods, we propose a
  linear optimization problem that computes the proposed detection
  measure and a semidefinite program that produces the proposed
  reachable set. %
\end{abstract}

\vspace*{-1ex}
\section{Introduction}
Cyber-Physical Systems (CPS) are physical processes that are tightly
integrated with computation and communication systems for monitoring
and control. Examples include critical infrastructure, such as
transportation networks, energy distribution systems, and the
Internet. These systems are usually complex, large-scale and
insufficiently supervised, making them vulnerable to
attacks~\cite{%
  AC-SA-BS-AG-AP-SS:09,FP-FD-FB:12a}. %
A significant literature has studied various \textit{denial of
  service}~\cite{SA-AC-SSS:09},~\textit{false
  data-injection}~\cite{FM-QZ-MP-GJP:16,CB-FP-VG:17},~\textit{replay}~\cite{YM-BS:09,MZ-SM:TAC11}, %
\textit{sensor, and integrity}
attacks~\cite{JM-DU-HS-KJ:18,YM-BS:15,SM-YS-NK-SD-PT:16,
  CM-NV-JR:17}. %
The majority of these works study attack-detection problems in a
control-theoretical framework. This approach essentially employs
detectors to identify abnormal behaviors by comparing estimation and
measurements under some predefined metrics. However, attacks could be
\textit{stealthy}, and exploit knowledge of the system structure,
uncertainty and noise information to inflict significant damage on the
physical system \reviseone{while avoiding} detection. This motivates
the characterization of the impact of stealthy attacks via
e.g.~reachability set
analysis~\cite{YM-BS:15,CB-FP-VG:15,CM-IS-JR-DN:18}.
To ensure computational tractability, these works assume either
Gaussian or bounded system noise. \reviseone{However, these assumptions
fall short in modeling all natural disturbances that
  can affect a system. Such systems would be vulnerable to stealthy
  attacks that disguise themselves via an intentionally selected, unbounded and non-Gaussian distribution.
  When designing detectors, an added difficulty is in
  obtaining tractable computations that can handle these more general
  distributions.}  More recently, novel measure of concentration has
opened the way for \reviseone{online tractable and} %
robust attack detection with probability guarantees \reviseone{under uncertainty}.
A first attempt in this direction
is~\cite{VR-NH-JR-TS:19}, where a Chebyshev inequality is used to
design a detector and, \reviseone{an assessment of} the impact of
stealthy attacks \reviseone{is given}, under the assumption that
system noises were bounded.
   \reviseone{With the aim of
  obtaining a less conservative detection mechanism,} we %
leverage \reviseone{an alternative} measure-concentration result
\reviseone{via Wasserstein metric. This metric
 is built from data gathered on the system, and can provide concentration result significantly sharper than that via Chebyshev inequality.} %
In particular, we %
address the following question for linear CPSs: \textit{How to design an
   \reviseone{online} attack-detection mechanism that is robust to
   \reviseone{light-tailed distributions of} system noise while
   remaining sensitive to attacks \reviseone{and limiting} the
   impact of the stealthy attack?}

 \reviseone{To answer the question}, we consider a sensor-attack
 detection problem on a linear dynamical system. The linear system
 models a remotely-observed process that is subject to an additive
 noise described by an unknown, not-necessarily bounded, light-tailed
 distribution. To identify abnormal behavior, we \reviseone{employ a
   steady-state Kalman filter} as well as a benchmark distribution of
 an offline sequence corresponding to the normal system
 operation. \reviseone{In this framework, we propose a novel detection
   mechanism that achieves online and robust attack detection of
   stealthy attacks in high confidence.}

 \textit{Statement of Contributions:}
 \reviseone{%
 1) We propose a novel}
 detection measure, \reviseone{which} employs the Wasserstein distance
 between the benchmark and a distribution of the residual sequence
 obtained online. \reviseone{2) We propose a novel}
 threshold-detection mechanism, \reviseone{which} exploits
 measure-of-concentration results to guarantee the robust detection of
 an attack with high confidence \reviseone{using a finite set of
   data}, and which further enables the robust tuning of the false
 alarm rate. The proposed detector can effectively identify real-time
 attacks when its behavior differs from that of the system
 noise. \reviseone{In addition, the detector can handle systems noises
   that are not necessarily distributed as Gaussian.} \reviseone{3) We
   propose a quantifiable, probabilistic state-reachable set, which
   reveals} the impact of \reviseone{the} stealthy attacker and system
 noise \reviseone{in high probability.}
 \reviseone{4)} To implement and analyze the proposed mechanism, we
 formulate a linear optimization problem and a semidefinite problem
 for the computation of the detection measure as well as the reachable
 set, respectively. 
We illustrate our methods in a two-dimensional linear system with irregular noise distributions and stealthy sensor attacks.

\vspace*{-1ex}
\section{CPS and Its Normal Operation}\label{sec:ProbStat}
This section presents cyber-physical system (CPS) model, and underlying
assumptions on the system and attacks. %
{
\begin{figure}[!tbp]%
\centering
\includegraphics[width=0.4\textwidth]{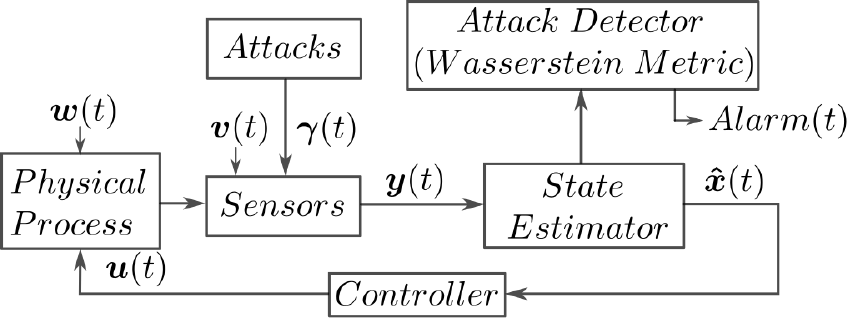}%
\caption{\small {Cyber-Physical System Diagram.}}%
\label{fig:cps}%
\end{figure}
}

A remotely-observed, cyber-physical system subject to
sensor-measurement attacks,
{as in Fig.~\ref{fig:cps}, }
is described as a discrete-time, stochastic, linear, and
time-invariant system
\begin{equation}
\begin{aligned}
  \vect{x}(t+1) &= \; A \vect{x}(t) + B \vect{u}(t) + \vect{w}(t) , \\
  \vect{y}(t)  \;  &= \; C \vect{x}(t)+ \vect{v}(t) + \vect{\gamma}(t),
\end{aligned} \label{eq:sys}
\end{equation}
where $\vect{x}(t) \in \real^n$, $\vect{u}(t) \in \real^m$ and
$\vect{y}(t) \in \real^p$ denote the system state, input and output at
time $t \in \mathbb{N}$, respectively. The state matrix $A$, input
matrix $B$ and output matrix $C$ are assumed to be known in
advance. In particular, we assume that %
the pair $(A,B)$ is stabilizable, and $(A,C)$ is
detectable.
The process noise $\vect{w}(t) \in \real^n$ and output noise
$\vect{v}(t) \in \real^p$ are independent zero-mean random vectors. We
assume that each $\vect{w}(t)$ and $\vect{v}(t)$ are independent and
identically distributed (i.i.d.) over time. %
We denote their (unknown, not-necessarily equal) distributions by
$\prob_{\vect{w}}$ and $\prob_{\vect{v}}$, respectively. In addition,
we assume that $\prob_{\vect{w}}$ and $\prob_{\vect{v}}$ are
light-tailed\footnote{\label{note:light-tail} For a random vector
  $\vect{w}$ such that %
  $\vect{w} \sim \prob$, we say $\prob$ is $q$-light-tailed,
  $q=1,2,\ldots$,
if %
$c:=\mathbb{E}_{\prob}[\exp( b\Norm{\vect{w}}^a)] < \infty$ for some $a>q$ and $b>0$.
All examples listed have a moment generating function, so their
exponential moment can be constructed for at least $q=1$.},
excluding scenarios of systems operating under extreme events, or
subject to large delays. In fact, Gaussian, Sub-Gaussian, Exponential
distributions, and any distribution with a compact support set are
admissible. This class of distributions is sufficient to characterize
the uncertainty or noise of many practical problems.

An additive sensor-measurement attack is implemented via
$\vect{\gamma}(t) \in \real^p$ in~\eqref{eq:sys}, on which we assume
the following
\begin{assumption}[Attack model] \label{assump:attack} It holds that
  1) $\vect{\gamma}(t) = \vect{0}$ whenever there is no attack;
  2) the attacker can modulate any component of $\vect{\gamma}(t)$ at
  any time; 3) the attacker has unlimited computational resources and
  access to system information, e.g., $A$, $B$, $C$, $\vect{u}$,
  $\prob_{\vect{w}}$ and $\prob_{\vect{v}}$ to decide on
  $\vect{\gamma}(t)$, $t \in \mathbb{N}$. %
\end{assumption}

\vspace*{-2ex}
\subsection{Normal System Operation}
\label{subsec:noattack}
In what follows, we introduce the state observer that enables
prediction in the absence of attacks (when
$\vect{\gamma}(t)=\vect{0}$). Since the distribution of system noise is
unknown, we identify a benchmark distribution that can be used to
characterize this unknown distribution with high confidence.

To predict (estimate) the system behavior, we leverage the system
information $(A, B, C)$ and employ a \reviseone{Kalman filter}%
\begin{equation*}
\begin{aligned}
  \hat{\vect{x}}(t+1) &= \; A \hat{\vect{x}}(t) + B \vect{u}(t) + L(t) \left( \vect{y}(t) - \hat{\vect{y}}(t)   \right) , \\
  \hat{\vect{y}}(t) \; &= \; C \hat{\vect{x}}(t),
\end{aligned} \label{eq:observer}
\end{equation*}
where $\hat{\vect{x}}(t)$ is the state estimate and $L(t) \equiv L$ is
the \reviseone{steady-state Kalman} %
gain matrix. As the pair $(A,C)$ is detectable, %
\reviseone{the} gain $L$ can be selected in such a way that the
eigenvalues of %
$A-LC$ are inside the unit circle. This
ensures asymptotic tracking of the state in expectation; that is, the
expectation of the estimation error $\vect{e}(t):=\vect{x}(t)-
\hat{\vect{x}}(t)$ satisfies
\begin{equation*}
  \mathbb{E} [\vect{e}(t)]  \rightarrow 0 \textrm{ as } t \rightarrow \infty, \; \textrm{ for any }  \vect{x}(0), \; \hat{\vect{x}}(0).
\end{equation*}
The further selection of eigenvalues of $A-LC$ and the structure of
$L$ usually depends on additional control objectives such as noise
attenuation requirements. In this paper, we additionally consider the
estimated state feedback
\begin{equation*}
\vect{u}(t)=K  \hat{\vect{x}}(t),
\end{equation*}
where $K$ is so that the following augmented system is
stable\footnote{\reviseone{System~\eqref{eq:augment} is
    input-to-state stable in probability (ISSp) relative to any
    compact set $\mathcal{A}$ which contains the origin, if we select
    $K$ such that eigenvalues of the matrix $A+BK$ are inside the unit
    circle, see e.g.~\cite{ART-JPH-AS:14}. }%
 }
\begin{equation} \label{eq:augment}
  \vect{\xi}(t+1)= F \vect{\xi}(t) + G \vect{\sigma} (t),
\end{equation}
where $\vect{\xi}(t):= \trans{ \left( {\vect{x}}(t) , \vect{e}(t)
  \right) }$, $\vect{\sigma}(t):= \trans{\left( \vect{w}(t) ,
    \vect{v}(t) +\vect{\gamma}(t) \right) }$,
\begin{equation*}
  F= \begin{bmatrix}
   A+BK & -BK \\ 0 & A-LC
 \end{bmatrix},
G= \begin{bmatrix}
 I & 0 \\ I & -L
\end{bmatrix}  \textrm{ and some } \vect{\xi}(0).
\end{equation*}
  \begin{remark}[Selection of $L$ and $K$] {\rm In general, the
      selection of the matrices $L$ and $K$ for the
      system~\eqref{eq:sys} is a nontrivial task, especially when
      certain performance criteria are to be satified, such as fast
      system response, energy conservation, or noise
      minimization. However, there are a few scenarios in which the
      \textit{Separation Principle} can be invoked for a tractable
      design of $L$ and $K$. For example, 1) when there is no system
      noise, matrices $L$ and $K$ can be designed separately, such
      that each $A+BK$ and $A-LC$ have all eigenvalues contained
      inside the unit circle, respectively. 2) when noise are
      Gaussian, the gain matrices $L$ and $K$ can be designed to
      minimize the steady-state covariance matrix and control effort,
      via a separated design of a Kalman filter (as an observer) and a
      linear-quadratic regulator (as a controller). The resulting
      design is referred to as a Linear-Quadratic-Gaussian (LQG)
      control~\cite{MA:71}.
  }\end{remark}

Consider the system initially operates normally with the proper
selection of $L$ and $K$, and assume that the augmented
system~\eqref{eq:augment} is in steady state, i.e.,
$\mathbb{E}[\vect{\xi}(t)]=\vect{0}$. %
In order to design an attack detector later, we need a
characterization of the distribution of the \textit{residue}
$\vect{r}(t) \in \real^p$, which measures the difference between what
we measure and what we expect to receive, as follows
\begin{equation*} \label{eq:r}
  \begin{aligned}
    \vect{r}(t) :=& \; \vect{y}(t) - \hat{\vect{y}}(t) = C\vect{e}(t)
    +\vect{v}(t) +\vect{\gamma}(t).
  \end{aligned}
\end{equation*}
When there is no attack, it can be verified that the random vector
$\vect{r}(t)$ is zero-mean, and light-tailed\footnote{This can be
  checked from the definition %
  in the
  footnote~\ref{note:light-tail}, and the fact that $\vect{r}(t)$ is a
  linear combination of zero-mean $q$-light-tailed distributions.},
and we denote its unknown distribution by $\prob_{\vect{r}}$. We
assume that a finite but large number $N$ of i.i.d.~samples of
$\prob_{\vect{r}}$ are accessible, and acquired by collecting
$\vect{r}(t)$ for a sufficiently large time. We call these
i.i.d.~samples a \textit{benchmark data set}, $\Xi_{\Batch}:=\{
\vect{r}^{(i)} \}_{i=1}^{N}$, and construct the resulting empirical
distribution $\prob_{\vect{r}, \Batch}$ by
\begin{equation*}
  \prob_{\vect{r}, \Batch}:= \frac{1}{N}\sum_{i=1}^{N} \delta_{ \{ \vect{r}^{(i)} \} },
\end{equation*}
where the operator $\delta$ is the mass function and the subscript
$\Batch$ indicates that $\prob_{\vect{r}, \Batch}$ is the benchmark to
approximate the unknown $\prob_{\vect{r}}$. We can claim that, the
benchmark $\prob_{\vect{r}, \Batch}$ provides a good sense of the
effect of the noise on the system~\eqref{eq:augment} via the following
measure concentration result
\begin{theorem}[Measure Concentration~{\cite[Application of
    Theorem~2]{NF-AG:15}}] %
  If $\prob_{\vect{r}}$ is a $q$-light-tailed distribution for some $q
  \geq 1$, then for a given $\beta \in (0,1)$, the following holds
\begin{equation*}
  {\Prob} \left( d_{W,q}(\prob_{\vect{r}, \Batch}, \prob_{\vect{r}} ) \leq \epsilon_{\Batch} \right)
  \geq 1- \beta,
  \label{eq:benchconcentrate}
\end{equation*}
where $\Prob$ denotes the Probability, %
$d_{W,q}$ denotes the $q$-Wasserstein
metric\footnote{\label{note:Wasserstein} Let
  ${\mathcal{M}_q}(\mathcal{Z})$ denote the space of all
  $q$-light-tailed probability distributions supported on $\mathcal{Z}
  \subset \real^p$. Then for any two distributions $\mathbb{Q}_1$,
  $\mathbb{Q}_2 \in \mathcal{M}_q(\mathcal{Z})$, the $q$-Wasserstein
  metric~\cite{FS:15} $\map{d_{W,q}}{\mathcal{M}_q(\mathcal{Z})
    \times \mathcal{M}_q(\mathcal{Z})}{\realnonnegative}$ is defined
  by
  \begin{equation*}
  d_{W,q}(\mathbb{Q}_1,\mathbb{Q}_2)
  := (  \min\limits_{\Pi} \int_{\mathcal{Z} \times \mathcal{Z}} \ell^q(\xi_1 , \xi_2) \Pi(d \xi_1, d \xi_2) )^{1/q},
  \end{equation*}
where $\Pi$ is in a set of all the probability distributions on
$\mathcal{Z} \times \mathcal{Z}$ with marginals $\mathbb{Q}_1$ and
$\mathbb{Q}_2$. The cost $\ell(\xi_1 , \xi_2):=\Norm{\xi_1 -\xi_2}$ is
a norm on $\mathcal{Z}$. } in the probability space,
and the parameter
$\epsilon_{\Batch}$ is selected as
\begin{equation}
  \epsilon_{\Batch} := \left\{ {\begin{array}{*{10}{l}}
        \left( \frac{\log(c_1 \beta^{-1})}{c_2 N} \right)^{q/a},
         & \textrm{if} \; N < \frac{\log(c_1 \beta^{-1})}{c_2},  \\
            \bar{\epsilon}  ,
           & \textrm{if} \; N \geq \frac{\log(c_1 \beta^{-1})}{c_2}, %
\end{array}} \right.
\label{eq:epsirad}
\end{equation}
 for some constant\footnote{The parameter $a$ is determined as in the
  definition of $\prob_{\vect{r}}$ and the constants $c_1$, $c_2$
  depend on $q$, $m$, and $\prob_{\vect{r}}$ (via $a$, $b$, $c$). When
  information on $\prob_{\vect{r}}$ is absent, the parameters $a$,
  $c_1$ and $c_2$ can be determined in a data-driven fashion using
  sufficiently many samples of $\prob_{\vect{r}}$. See~\cite{NF-AG:15} for details.} $a>q$, $c_1$,
$c_2>0$, and $\bar{\epsilon}$ is such that
 \begin{equation*} \small
 \small \frac{\bar{\epsilon}}{ \log (2 + 1/\bar{\epsilon}) } =
  \left( \frac{\log(c_1 \beta^{-1})}{c_2 N} \right)^{1/2}, \textrm{ if } p =2q,
 \end{equation*}
or
 \begin{equation*}
\small  \bar{\epsilon}=
  \left( \frac{\log(c_1 \beta^{-1})}{c_2 N} \right)^{1/\max\{2,p/q \}}, \textrm{ if } p \ne 2q,
 \end{equation*}
where $p$ is the dimension of $\vect{r}$. \hfill $\square$
\label{thm:concentration}
\end{theorem}

Theorem~\ref{thm:concentration} provides a probabilistic bound
$\epsilon_{\Batch}$ on the $q$-Wasserstein distance between
$\prob_{\vect{r}, \Batch}$ and $\prob_{\vect{r}}$, with a confidence
at least $1-\beta$. %
The result indicates how to tune the parameter $\beta$ and the number
of benchmark samples $N$ that are needed to find a sufficiently good
approximation of $\prob_{\vect{r}}$, by means of $\prob_{\vect{r},
  \Batch}$. In this way, given a fixed $\epsilon$, we can increase our
confidence ($1-\beta$) on whether $\prob_{\vect{r}}$ and
$\prob_{\vect{r}, \Batch}$ are within distance $\epsilon$, by
increasing the number of samples. We assume that $\prob_{\vect{r},
  \Batch}$ has been determined in advance, selecting a very large
(unique) $N$ to ensure various very small bounds $\epsilon_{\Batch}$
associated with various $\beta$.
{Later, we discuss how
  the parameter $\beta$ can be interpreted as a \textit{false alarm
    rate} in the proposed attack detector. The resulting
  $\prob_{\vect{r}, \Batch}$, with a tunnable false alarm rate
  (depending on $\beta$), will allow us to design a detection
  procedure which is robust to the system noise.}

\section{Threshold-based robust detection of attacks, and
  stealthiness} \label{sec:attack_detection}

This section presents our online detection procedure, and a
threshold-based detector with high-confidence performance guarantees. Then, we propose a tractable computation of the detection
measure used for online detection. We finish the section by
introducing a class of stealthy attacks.

\noindent \textbf{Online Detection Procedure (ODP):} At each time $t\ge T$, we
construct a $T$-step detector distribution %
\begin{equation*}
  \prob_{\vect{r}, \Detect}:= \frac{1}{T}\sum_{j=0}^{T-1} \delta_{ \{ \vect{r}(t-j) \} },
\end{equation*}
where $\vect{r}(t-j)$ is the residue data collected independently at
time $t-j$, for $j \in \untilinterval{0}{T-1}$. Then with a given $q$
and a threshold $\alpha >0$, we consider the \textit{detection
  measure}
\begin{equation} \label{eq:z}
  z(t):= d_{W,q}(\prob_{\vect{r}, \Batch}, \prob_{\vect{r}, \Detect} ),
\end{equation}
and the \textit{attack detector}
\begin{equation} \label{eq:detector}
  \begin{cases}
    z(t) \leq \alpha, & \textrm{ no alarm at }t: \Alarm(t)=0,
    \\
      z(t) > \alpha, & \textrm{ alarm at }t: \Alarm(t)=1  ,
  \end{cases}
\end{equation}
with $\Alarm(t)$ the sequence of alarms generated online based on
the previous threshold.

The distribution $\prob_{\vect{r}, \Detect}$ uses a small number $T$
of samples to ensure the online computational tractability of $z(t)$,
so $\prob_{\vect{r}, \Detect}$ is highly dependent on instantaneous
samples. Thus, $\prob_{\vect{r}, \Detect}$ may significantly deviate
from the true $\prob_{\vect{r}}$, and from $\prob_{\vect{r},
  \Batch}$. Therefore, even if there is no
attack, %
the attack detector is expected to generate false alarms due to the
system noise as well as an improper selection of the threshold
$\alpha$. In the following, we discuss how to select an $\alpha$ that
is robust to the system noise and which results in a desired false
alarm rate. Note that the value $\alpha$ should be small to be able to
distinguish attacks from noise, as discussed later.
\begin{lemma}[Selection of $\alpha$ for Robust Detectors]
  Given parameters $N$, $T$, $q$, $\beta$, and a desired false alarm
  rate $\Delta>\beta$ at time $t$, if we select the threshold $\alpha$
  as
\begin{equation*}
  \alpha:= \epsilon_{\Batch} +\epsilon_{\Detect},
\end{equation*}
where $\epsilon_{\Batch}$ is chosen as in~\eqref{eq:epsirad} and
$\epsilon_{\Detect}$ is selected following the
$\epsilon_{\Batch}$-formula~\eqref{eq:epsirad}, but with $T$ and
$\frac{\Delta-\beta}{1-\beta}$ in place of $N$ and $\beta$,
respectively. Then, the detection measure~\eqref{eq:z} satisfies
\begin{equation*}
  {\Prob} \left( z(t) \leq \alpha \right)
  \geq 1- \Delta,
\label{eq:robust2noise}
\end{equation*}
for any zero-mean $q$-light-tailed underlying distribution
$\prob_{\vect{r}}$.
\label{lemma:robust2noise}
\end{lemma} %
  Due to space limit, please see ArXiv version~\cite{DL-SM:20-extended} for proofs.

\begin{proof}
  The proof leverages the  triangular inequality
\begin{equation*}
z(t) \leq d_{W,q}(\prob_{\vect{r}, \Batch}, \prob_{\vect{r}} ) + d_{W,q}(\prob_{\vect{r}, \Detect}, \prob_{\vect{r}} ),
\end{equation*}
the measure concentration result for each $d_{W,q}$ term, and that
samples of $\prob_{\vect{r}, \Batch}$ and $\prob_{\vect{r}, \Detect}$
are collected independently.
\end{proof}
Lemma~\ref{lemma:robust2noise} ensures that the false alarm rate is no
higher than $\Delta$ when there is no
attack, %
i.e.,
\begin{equation*}
  \Prob ( \Alarm(t)=1 \; | \; \noattack) \leq \Delta, \quad \forall \; t.
\end{equation*}
Note that the rate $\Delta$ can be selected by properly choosing the
threshold $\alpha$. Intuitively, if we fix all the other parameters,
then the smaller the rate $\Delta$, the larger the threshold
$\alpha$. Also, large values of $N$, $T$, $1-\beta$ contribute to
small $\alpha$.
\begin{remark}[Comparison with $\chi^2$-detector]
Consider an alternative detection measure
\begin{equation*}
  z_{\chi}(t):= \trans{\vect{r}(t)} \Sigma^{-1} \vect{r}(t),
\end{equation*}
where $\Sigma$ is the constant covariance matrix of the residue $\vect{r}(t)$ under normal system operation. In particular, if $\vect{r}$ is Gaussian, the detection measure $z_{\chi}(t)$ is $\chi^2$-distributed and referred to as $\chi^2$ detection measure with $p$ degree of freedom. The detector threshold $\alpha$ is selected via look-up tables of $\chi^2$ distribution, given the desired false alarm rate $\Delta$.
To compare $z(t)$ with $z_{\chi}(t)$, we leverage the assumption that $\vect{r}$ is Gaussian with the given covariance $\Sigma$.
This gives explicitly the expression of $z(t)$ the following
\begin{equation*}
  z(t):=  \left( \mathbb{E}_{\xi \sim \mathcal{N}(\vect{r}(t), \Sigma ) } [ \Norm{ \xi }^q  ]  \right)^{-1/q}.
\end{equation*}
By selecting $q=2$, we have
\begin{equation*}
  z(t):= \left( \trans{\vect{r}(t)} \vect{r}(t) + Tr(\Sigma)\right)^{-1/2}.
\end{equation*}
Note that, the measure-of-concentration result in Theorem~\ref{thm:concentration} is sharp when $\vect{r}$ is Gaussian, which in fact results in the threshold $\alpha$ as tight as that derived for $\chi^2$-detector.
\end{remark}

\noindent \textbf{Computation of detection measure:} To achieve a
tractable computation of the detection measure $z(t)$, we leverage the
definition of the Wasserstein distance (see
footnote~\ref{note:Wasserstein}) and the fact that both
$\prob_{\vect{r}, \Batch}$ and $\prob_{\vect{r}, \Detect}$ are
discrete. The solution is given as a linear program.

The Wasserstein distance $d_{W,q}(\prob_{\vect{r}, \Batch},\prob_{\vect{r}, \Detect})$,
originally solving the \textit{Kantorovich optimal transport
  problem}~\cite{FS:15}, can be interpreted as the minimal work needed
to move a mass of points described via a probability distribution $\prob_{\vect{r},
  \Batch}(\vect{r})$, on the space $\mathcal{Z} \subset \real^p$, to a
mass of points described by the probability distribution
$\prob_{\vect{r}, \Detect}(\vect{r})$ on the same space, with some
transportation cost $\ell$.  The minimization that defines $d_{W,q}$
is taken over the space of all the joint distributions $\Pi$ on
$\mathcal{Z} \times \mathcal{Z}$ whose marginals are $\prob_{\vect{r},
  \Batch}$ and $\prob_{\vect{r}, \Detect}$, respectively.
Assuming that both $\prob_{\vect{r}, \Batch}$ and $\prob_{\vect{r},
  \Detect}$ are discrete, we can equivalently characterize the joint
distribution $\Pi$ as a discrete mass \textit{optimal transportation
  plan}~\cite{FS:15}. To do this,
let us consider two sets $\mathcal{N}:=\untilinterval{1}{N}$ and
$\mathcal{T}:=\untilinterval{0}{T-1}$. Then, $\Pi$ can be
parameterized (by $\lambda$) as follows
\begin{align}
  \Pi_{\lambda}&(\xi_1,\xi_2):= \sum_{i\in\mathcal{N}} \sum_{j\in\mathcal{T}} \lambda_{ij} \delta_{ \{ \vect{r}^{(i)} \} }(\xi_1) \delta_{ \{ \vect{r}(t-j) \} } (\xi_2), \nonumber \\
  \st & \hspace{0cm} \sum_{i\in\mathcal{N}}\lambda_{ij} = \frac{1}{T}, \; \forall \; j \in \mathcal{T}, \quad \sum_{j\in\mathcal{T}} \lambda_{ij}=\frac{1}{N}, \; \forall \; i \in \mathcal{N}, \label{eq:lambda1} \\
  & \hspace{0.3cm} \lambda_{ij} \geq 0, \; \forall \; i \in
  \mathcal{N}, \; j \in \mathcal{T}. \label{eq:lambda2}
\end{align}
Note that this characterizes all the joint distributions with
marginals $\prob_{\vect{r}, \Batch}$ and $\prob_{\vect{r}, \Detect}$,
where $\lambda$ is the allocation of the mass from $\prob_{\vect{r},
  \Batch}$ to $\prob_{\vect{r}, \Detect}$. Then, the proposed
detection measure $z(t)$ in~\eqref{eq:z} reduces to the following
\begin{equation}
  \begin{aligned}
    (z(t) )^q:=   \min\limits_{\lambda}  & \;  \sum_{i\in \mathcal{N}} \sum_{j\in \mathcal{T}} \lambda_{ij} \Norm{\vect{r}^{(i)} - \vect{r}(t-j)}^q , \\
    \st \; & \; \eqref{eq:lambda1}, \; \eqref{eq:lambda2},
  \end{aligned} \label{eq:P} \tag{P}
\end{equation}
which is a linear program over a compact polyhedral set. Therefore,
the solution exists and~\eqref{eq:P} can be solved to global optimal
in polynomial time by e.g., a CPLEX  solver.

\subsection{Detection and Stealthiness of
  Attacks} \label{subsec:attack} Following from the previous
discussion, we now introduce a False Alarm Quantification Problem,
then specialize it to the Attack Detection Problem considered in this
paper. In particular, we analyze the sensitivity of the proposed
attack detector method for the cyber-physical system under attacks.

\noindent \textbf{Problem 1. (False Alarm Quantification Problem):}
Given the augmented system~\eqref{eq:augment},
the online detection procedure in
Section~\ref{sec:attack_detection}, and the attacker type described in
Assumption~\ref{assump:attack}, compute the false alarm rate
\begin{equation*}
\begin{aligned}
&  \Prob(\textrm{false alarm at }t):= \\
& \hspace{1cm} \Prob ( \Alarm(t)=1 \; | \; \noattack) \Prob (  \noattack) \\
& \hspace{1.9cm} + \Prob ( \Alarm(t)=0 \; | \; \attack) \Prob (\attack).
\end{aligned}
\end{equation*}
Problem~1, on the computation of the false alarm probability, requires
prior information of the attack probability $\Prob
(\attack)$. In this work, we are interested in
\textit{stealthy attacks}, i.e., attacks that cannot be effectively
detected by~\eqref{eq:detector}. We are led to the following problem.

\noindent \textbf{Problem 2. (Attack Detection Problem):} Given the
setting of Problem 1, provide conditions that characterize stealthy
attacks, i.e., attacks that contribute to $\Prob ( \Alarm(t)=0
\; | \; \attack)$, and quantify their impact on the system.

To remain undetected, the attacker must select $\vect{\gamma}(t)$ such
that $z(t)$ is limited to within the threshold
$\alpha$. %
To quantify the effects of these attacks, let us consider an attacker
sequence backward in time with length $T$. At time $t$, denote the
arbitrary injected attacker sequence by $\vect{\gamma}(t-j) \in
\real^p$, $j \in \untilinterval{0}{T-1}$ (if $t-j<0$, assume
$\vect{\gamma}(t-j)=0$).
This sequence, together with~\eqref{eq:augment}, results in a detection sequence
$\{\vect{r}(t-j)\}_j$
that can be used to construct detector
distribution $\prob_{\vect{r}, \Detect}$ and detection measure
$z(t)$. For simplicity and w.l.o.g., let us assume that $\vect{\gamma}(t)$ is in the following form
\begin{equation} \label{eq:gammabar} \vect{\gamma}(t):= - C\vect{e}(t)
  - \vect{v}(t) + \bar{\gamma}(t),
\end{equation}
where $\bar{\gamma}(t) \in \real^p$ is any vector selected by the attacker.
\footnote{
\reviseone{Note that, when there is no attack at $t$, we have $\vect{\gamma}(t) = \vect{0}$, resulting in $\bar{\gamma}(t)= C\vect{e}(t) + \vect{v}(t)$. Similar techniques appearred in, e.g.,~\cite{YM-BS:15,CM-JR:16a}.}}
We characterize the scenarios that can occur, providing a first,
partial answer to Problem~2. Then, we will come back to characterizing
the impact of stealthy attacks in
Section~\ref{sec:stealth}. %

\begin{definition}[Attack Detection Characterization]
  Assume~\eqref{eq:augment} is subject to attack, i.e.,
  $\vect{\gamma}(t) \ne \vect{0}$
    for some $t \geq 0$.
\begin{itemize}
\item If $z(t) \leq \alpha$, $\forall \; t \geq 0$, then the attack is
  stealthy with probability one, i.e., $\Prob ( \Alarm(t)=0 \; |
   \attack) =1$.
\item If $z(t) \leq \alpha$, $\forall t \leq M$, then the attack is
  $M$-step stealthy.
\item If $z(t) > \alpha$, $\forall t \geq 0$, then the attack is
  active with probability one, i.e., $\Prob ( \Alarm(t)=0 \; |
   \attack) =0$.
\end{itemize}
\label{lemma:asattack}
\end{definition}
\begin{lemma}[Stealthy Attacks Leveraging System Noise]
  Assume~\eqref{eq:augment} is subject to attack in
  form~\eqref{eq:gammabar}, where $\bar{\gamma}(t)$ is stochastic and
  distributed as $\prob_{\bar{\gamma}}$. If $\prob_{\bar{\gamma}}$ is
  selected such that $d_{W,q}(\prob_{\bar{\gamma}} , \prob_{\vect{r},
    \Batch} ) \leq \epsilon_{\Batch}$, then the attacker is stealthy
  with (high) probability at least $\frac{1- \Delta}{1-\beta}$, i.e.,
  $\Prob ( \Alarm(t)=0 \; | \; \attack ) \geq \frac{1-
    \Delta}{1-\beta}$.
\end{lemma}
\begin{proof}
  Assume~\eqref{eq:augment} is under attack. Then, we prove the
  statement leveraging the measure concentration
\begin{equation*}
  {\Prob} \left( d_{W,q}(\prob_{\bar{\gamma}}, \prob_{\vect{r}, \Detect}) \leq \epsilon_{\Detect} \right)
  \geq 1- \frac{\Delta-\beta}{1-\beta},
\end{equation*}
which holds as $\prob_{\vect{r}, \Detect}$ is constructed using
samples of $\prob_{\bar{\gamma}}$,
and the triangular inequality
\begin{equation*}
  z(t) \leq d_{W,q}(\prob_{\vect{r}, \Batch}, \prob_{\bar{\gamma}} ) + d_{W,q}(\prob_{\vect{r}, \Detect}, \prob_{\bar{\gamma}} ),
\end{equation*}
resulting in $z(t) \leq \alpha$ with probability at least $\frac{1-
  \Delta}{1-\beta}$. %
\end{proof}
\reviseone{%
Note that $\alpha>\epsilon_{\Batch}$, which allows %
the attacker to select $\prob_{\bar{\gamma}}$ with $\epsilon_{\Batch} < d_{W,q}(\prob_{\bar{\gamma}} , \prob_{\vect{r},\Batch} ) \leq \alpha$.
However, the probability of being stealthy can be indefinitely low, with the range $[0,\frac{1- \Delta}{1-\beta}]$.
}

\section{Stealthy Attack Analysis} \label{sec:stealth} In this
section, we address the second question in Problem~2 via reachable-set
analysis. In particular,
we achieve an outer-approximation of the finite-step probabilistic
reachable set, quantifying the effect of the stealthy attacks and the
system noise in probability.

Consider an attack sequence $\vect{\gamma}(t)$ as
in~\eqref{eq:gammabar}, where $\bar{\gamma}(t) \in \real^p$ is any
vector such that the attack remains stealthy. That is,
$\bar{\gamma}(t)$ results in the detected distribution
$\prob_{\vect{r}, \Detect}$, which is close to $\prob_{\vect{r},
  \Batch}$ as prescribed by $\alpha$. This results in the following
representation of the system~\eqref{eq:augment}
\begin{equation} \label{eq:attack} \small
  \vect{\xi}(t+1)=
  \underbrace{\begin{bmatrix}
   A+BK & -BK \\ 0 & A
 \end{bmatrix}}_{H} \vect{\xi}(t) + \underbrace{\begin{bmatrix}
  I & 0 \\ I & -L
\end{bmatrix}}_{G} \begin{bmatrix}
   \vect{w}(t) \\
   \bar{\gamma}(t)
\end{bmatrix}.
\end{equation}
\reviseone{
 We provide in the following remark an intuition of how restrictive the stealthy attacker's action $\bar{\gamma}(t)$ has to be.
\begin{remark}[Constant attacks]
  {\rm Consider a constant offset attack $\bar{\gamma}(t):= \gamma_0$ for
    some $\gamma_0 \in \real^p$, $\forall t$.  Then
    by~\eqref{eq:P}, %
  \begin{equation*}
    z(t)=N^{1-1/q}\Norm{\gamma_0 - \textrm{C}(\Xi_{\Batch})}, \quad \textrm{C}(\Xi_{\Batch}):= \frac{1}{N} \sum_{i\in \mathcal{N}} \vect{r}^{(i)}.
  \end{equation*}
  To ensure stealth, we require $z(t) \leq
  \alpha$, this then limits the
  selection of $\gamma_0$ in a ball centered at
  $\textrm{C}(\Xi_{\Batch})$ with radius $\alpha/N^{1-1/q}$. Note that the radius can be
  arbitrarily small by choosing the benchmark size $N$ large.  }
\end{remark}
}

To quantify the reachable set of the system under attacks, prior
information on the process noise $\vect{w}(t)$ is needed.
To characterize $\vect{w}(t)$, let us assume that,
similarly to the benchmark $\prob_{\vect{r}, \Batch}$, we are able to
construct a noise benchmark distribution, denoted by $\prob_{\vect{w},
  \Batch}$. As before, %
\begin{equation*}
  {\Prob} \left( d_{W,q}(\prob_{\vect{w}, \Batch}, \prob_{\vect{w}} ) \leq \epsilon_{\vect{w},\Batch} \right)
  \geq 1- \beta,
\label{eq:wconcentrate}
\end{equation*}
for some $\epsilon_{\vect{w},\Batch}$. Given certain time, we are
interested in where, with high probability, the state of the system
can evolve from some $\vect{\xi}_0$. To do this, we consider the \textit{$M$-step
  probabilistic reachable set} defined as follows
\begin{equation*} \small
  \begin{aligned}
\mathcal{R}_{\vect{x},M}:= \left\{ \hspace{-0.2ex} \vect{x}(M)\in \real^n \; \rule[-0.4cm]{0.02cm}{1cm} \;
\begin{array}{l}
   \textrm{system } \eqref{eq:attack} \textrm{ with } \vect{\xi}(0)= \vect{\xi}_0, \\ %
  \exists \; \prob_{\vect{w}}  \ni d_{W,q}(\prob_{\vect{w}}, \prob_{\vect{w}, \Batch} ) \leq \epsilon_{\vect{w},\Batch}, \\
  \exists \; \prob_{\bar{\gamma}}  \ni d_{W,q}(\prob_{\bar{\gamma}} , \prob_{\vect{r}, \Batch} )  \leq \alpha , \\
\end{array}  \right\},
  \end{aligned}
\end{equation*}
then the true system state $\vect{x}(t)$ at time $M$, $\vect{x}(M)$,
satisfies
\begin{equation*}
  \Prob \left( \vect{x}(M) \in \mathcal{R}_{\vect{x},M}  \right) \geq 1-\beta,
\end{equation*}
where %
$1-\beta$ accounts for the independent restriction of the unknown
distributions $\prob_{\vect{w}}$ to be ``close'' to its benchmark.

The exact computation of $\mathcal{R}_{\vect{x},M}$ is intractable due
to the unbounded support of the unknown distributions $\prob_{\vect{w}}$
and $\prob_{\bar{\gamma}}$, even if they are close to their
benchmark. To ensure a tractable approximation, we follow a two-step
procedure. First, we equivalently characterize the constraints on
$\prob$ by its \textit{probabilistic support set}. Then, we
outer-approximate the probabilistic support by ellipsoids, and then the
reachable set by an ellipsoidal bound.

\noindent \textbf{Step 1: (Probabilistic Support of
  $\prob_{\bar{\gamma}} \; \ni d_{W,q}(\prob_{\bar{\gamma}} ,
  \prob_{\vect{r}, \Batch} ) \leq \alpha$)} We achieve this by
leveraging 1) the \textit{Monge formulation}~\cite{FS:15} of optimal
transport, 2) the fact that $\prob_{\vect{r}, \Batch}$ is discrete, and 3) results from coverage
control~\cite{FB-JC-SM:09,JC:08-tac}.
W.l.o.g., let us assume $\prob_{\bar{\gamma}}$ is non-atomic (or continuous) and,
consider the
distribution $\prob_{\bar{\gamma}}$ and
$\prob_{\vect{r}, \Batch}$ supported on
$\mathcal{Z}\subset \real^p$. Let us denote
by %
$\maps{f}{\prob_{\bar{\gamma}}}{\prob_{\vect{r}, \Batch}}$ the
\textit{transport map} that assigns mass over $\mathcal{Z}$ from
$\prob_{\bar{\gamma}}$ to $\prob_{\vect{r}, \Batch}$. The Monge
formulation of optimal transport aims to find an optimal transport map
that minimizes the transportation cost $\ell$ as follows
 \begin{equation*}
   d_{M,q}(\prob_{\bar{\gamma}}, \prob_{\vect{r}, \Batch})
   := \left(  \min\limits_{f} \int_{\mathcal{Z}} \ell^q(\xi , f(\xi) ) \prob_{\bar{\gamma}}(\xi) d \xi  \right)^{1/q}.
 \end{equation*}
 It is known that if an optimal transport map $f^{\star}$ exists, then
 the optimal transportation plan $\Pi^{\star}$ exists and the two
 problems $d_{M,q}$ and $d_{W,q}$ coincide\footnote{This is because
   the Kantorovich transport problem is the tightest relaxation of the
   Monge transport problem. See e.g.,~\cite{FS:15} for details.}. In
 our setting, for strongly convex $\ell^p$, and by the fact
 that $\prob_{\bar{\gamma}}$ is absolutely continuous, a
   unique optimal transport map can indeed be
   guaranteed\footnote{The Monge formulation is not always
   well-posed, i.e., there exists optimal transportation plans
   $\Pi^{\star}$ while transport map does not exist~\cite{FS:15}.},
 and therefore, $d_{M,q}=d_{W,q}$.

 Let us now consider any transport map $f$ of $d_{M,q}$, and define a
 partition of the support of $\prob_{\bar{\gamma}}$ by
\begin{equation*}
  W_i:=\setdef{\vect{r} \in \mathcal{Z}}{ f(\vect{r})=\vect{r}^{(i)}}, \; i \in \mathcal{N},
\end{equation*}
where $\vect{r}^{(i)}$ are samples in $\Xi_{\Batch}$, which generate $\prob_{\vect{r}, \Batch}$. %
Then, we equivalently
rewrite %
the value of the objective function defined in $d_{M,q}$, as
\begin{equation*} \small
  \begin{aligned}
    &  \mathcal{H}(\prob_{\bar{\gamma}}, W_1,\ldots,W_N ):= \sum\limits_{i=1}^N \int_{W_i} \ell^q(\xi , \vect{r}^{(i)} ) \prob_{\bar{\gamma}}(\xi) d \xi,  \\
  \end{aligned}
\end{equation*}
\begin{align} \small \label{eq:mass} \hspace{1cm} \st \; \int_{W_i}
  \prob_{\bar{\gamma}}(\xi) d \xi = \frac{1}{N}, \; \forall \; i \in
  \mathcal{N},
\end{align}
where the $\supscr{i}{th}$ constraints come from the fact that a
transport map $f$ should lead to consistent calculation of set volumes
under $\prob_{\vect{r}, \Batch}$ and $\prob_{\bar{\gamma}}$,
respectively. This results in the following equivalent problem to
$d_{M,q}$ as
\begin{equation}
  \begin{aligned}
    (d_{M,q}(\prob_{\bar{\gamma}}, \prob_{\vect{r}, \Batch}))^q:=   \min\limits_{W_i, i \in \mathcal{N}}  & \;  \mathcal{H}(\prob_{\bar{\gamma}}, W_1,\ldots,W_N ) , \\
    \st \; & \; \eqref{eq:mass}.
  \end{aligned} \label{eq:P1} \tag{P1}
\end{equation}
Given the distribution $\prob_{\bar{\gamma}}$, Problem~\eqref{eq:P1}
reduces to a load-balancing problem as in~\cite{JC:08-tac}. Let us
define the Generalized Voronoi Partition (GVP) of $\mathcal{Z}$
associated to the sample set $\Xi_{\Batch}$ and weight $\omega \in
\real^N$,
 for all $i \in \mathcal{N}$, as
\begin{equation*} \small
  \begin{aligned}
    & \mathcal{V}_i(\Xi_{\Batch}, \omega):=  \\
    & \hspace{0.8cm} \setdef{\xi \in \mathcal{Z}}{ \Norm{\xi -
        \vect{r}^{(i)}}^q - \omega_i \leq \Norm{\xi -
        \vect{r}^{(j)}}^q - \omega_j , \; \forall j \in \mathcal{N} }.
  \end{aligned}
\end{equation*}
It has been established that the optimal Partition of~\eqref{eq:P1} is the
  GVP~{\cite[Proposition V.1]{JC:08-tac}}. Further, the standard
  Voronoi partition,
  i.e., the GVP with equal weights $\bar{\omega}:= \vect{0}$,
  results in a lower
     bound%
  of~\eqref{eq:P1}, when constraints~\eqref{eq:mass} are
  removed~\cite{FB-JC-SM:09}, and therefore that of $d_{M,q}$. We
  denote this lower bound as
  $L(\prob_{\bar{\gamma}},\mathcal{V}(\Xi_{\Batch}) )$,
  and use this to quantify  a probabilistic
  support %
  of $\prob_{\bar{\gamma}}$. Let us consider the support set
\begin{equation*}
  \Omega(\Xi_{\Batch},\epsilon):= \cup_{i \in \mathcal{N}} \left( \mathcal{V}_i(\Xi_{\Batch}) \cap \mathbb{B}_{\epsilon}(\vect{r}^{(i)}) \right),
\end{equation*}
where %
{$\mathbb{B}_{\epsilon}(\vect{r}^{(i)}):=\setdef{\vect{r}\in
    \real^p}{ \Norm{\vect{r} - \vect{r}^{(i)}} \leq \epsilon}$.}
Then we have the following lemma.
\begin{lemma}[Probabilistic Support]
  Let $\epsilon >0$ and let $\prob_{\bar{\gamma}}$ be a distribution
  such that
  $L(\prob_{\bar{\gamma}},\mathcal{V}(\Xi_{\Batch})) \leq
  \epsilon^q$.
  Then, for any given $s>1$, at least $1-1/s^q$ portion
  of the mass of $\prob_{\bar{\gamma}}$ is supported on
  $\Omega(\Xi_{\Batch}, s\epsilon)$, i.e.,
  $\int_{\Omega(\Xi_{\Batch},s\epsilon)} \prob_{\bar{\gamma}}(\xi) d
  \xi \geq 1- 1/s^q$.
\end{lemma}
\begin{proof}
Suppose otherwise, i.e., $\int_{ \real^p \setminus \Omega(\Xi_{\Batch},s\epsilon)} \prob_{\bar{\gamma}}(\xi) d \xi = 1- \int_{\Omega(\Xi_{\Batch},s\epsilon)} \prob_{\bar{\gamma}}(\xi) d \xi > 1/s^q$. Then,
\begin{equation*} \small
  \begin{aligned}
 L(\prob_{\bar{\gamma}},\mathcal{V}(\Xi_{\Batch})) & \geq
    \int_{ \real^p \setminus \Omega(\Xi_{\Batch}, s\epsilon)} \Norm{\xi - \vect{r}^{(i)}}^q \prob_{\bar{\gamma}}(\xi) d \xi, \\
&  \geq s^q \epsilon^q \int_{ \real^p \setminus \Omega(\Xi_{\Batch},s\epsilon)} \prob_{\bar{\gamma}}(\xi) d \xi > \epsilon^q.
  \end{aligned}
\end{equation*}
\end{proof}
In this way, the support $\Omega(\Xi_{\Batch}, 2\alpha)$ contains at
least $1-1/2^q$ of the mass of all the distributions
$\prob_{\bar{\gamma}}$ such that $d_{W,q}(\prob_{\bar{\gamma}} ,
\prob_{\vect{r}, \Batch} ) \leq \alpha$. Equivalently, we have
\begin{equation*}
  \Prob \left( \bar{\gamma} \in  \Omega(\Xi_{\Batch}, 2\alpha) \right) \geq 1- 1/2^q,
\end{equation*}
where the random variable $\bar{\gamma}$ has a distribution $\prob_{\bar{\gamma}}$ such that $d_{W,q}(\prob_{\bar{\gamma}}, \prob_{\vect{r}, \Batch} ) \leq \alpha$. We characterize
$\prob_{\vect{w}}$ similarly.

Note that in practice, one can choose ball radius factor $s$ large
  in order to generate support which contains higher portion of the mass of unknown $\prob$. However, it comes at a cost on the approximation of
  the reachable set.

\noindent \textbf{Step 2: (Outer-approximation of
  $\mathcal{R}_{\vect{x},M}$)} Making use of the probabilistic
support, we can now obtain a finite-dimensional characterization
of the probabilistic reachable set, as follows
\begin{equation*} \small
  \begin{aligned}
\mathcal{R}_{\vect{x},M}:= \left\{ \hspace{-0.2ex} \vect{x}(M)\in \real^n \; \rule[-0.4cm]{0.02cm}{1cm} \;
\begin{array}{l}
  \textrm{system } \eqref{eq:attack}, \; \vect{\xi}(0)= \vect{\xi}_0, \\ %
  \vect{w} \in \Omega(\Xi_{\vect{w},\Batch}, 2\epsilon_{\vect{w},\Batch}), \\
  \bar{\gamma} \in  \Omega(\Xi_{\Batch}, 2\alpha)
\end{array}  \right\},
  \end{aligned}
\end{equation*}
and the true system state $\vect{x}(t)$ at time $M$, $\vect{x}(M)$,
satisfies
$  \Prob \left( \vect{x}(M) \in \mathcal{R}_{\vect{x},M}  \right) \geq (1-\beta)(1-1/2^q)^2.$
Many works focus on the tractable evolution of geometric
shapes; e.g.~\cite{YM-BS:15,CM-IS-JR-DN:18}. Here,
we follow~\cite{CM-IS-JR-DN:18} and propose outer
ellipsoidal bounds for $\mathcal{R}_{\vect{x},M}$. Let $Q_{\vect{w}}$
be the positive-definite shape matrix such that
$\Omega(\Xi_{\vect{w},\Batch},\epsilon_{\vect{w},\Batch}) \subset
\mathcal{E}_{\vect{w}}:=\setdef{\vect{w}}{ \trans{\vect{w}}
  Q_{\vect{w}}\vect{w} \leq 1 }$. Similarly, we denote
$Q_{\bar{\gamma}}$ and $\mathcal{E}_{\bar{\gamma}}$ for that of
$\bar{\gamma}$.  We now state the following lemma, that
applies~\cite[Proposition 1]{CM-IS-JR-DN:18} for our case.
\begin{lemma}[Outer bounds of
  $\mathcal{R}_{\vect{x},M}$]
  Given any $a \in [0,1)$, we claim $\mathcal{R}_{\vect{x},M} \subset
  \mathcal{E}(Q):=\setdef{\vect{x}\in\real^n}{
    \trans{\vect{\xi}}Q\vect{\xi} \leq a^{M}
    \trans{\vect{\xi}_0}Q\vect{\xi}_0 + \frac{(2-a)(1-a^{M})}{1-a}}$,
  with $Q$ satisfying
\begin{equation}
  \begin{aligned}
    Q \succ 0, \;\;
\begin{bmatrix}
  aQ & \trans{H}Q & \bf{0} \\
  QH & Q  & QG \\
  \bf{0} & \trans{G}Q & W
\end{bmatrix} \succeq 0, \\
\end{aligned} \label{eq:Q}
\end{equation}
where $H$, $G$ are that in~\eqref{eq:attack} and
\begin{equation}
\begin{aligned}
  W= \begin{bmatrix}
    (1-a_1)Q_{\vect{w}} & \bf{0} \\
    \bf{0} & (1-a_2) Q_{\bar{\gamma}}
  \end{bmatrix}, \\
  \textrm{for some } a_1+ a_2 \geq a, \; a_1, a_2 \in [0,1).
\end{aligned} \label{eq:W}
\end{equation}
\end{lemma}
\begin{figure*}[htp]
     \begin{minipage}[l]{0.65\columnwidth}
         \centering
         \includegraphics[width=1.05\textwidth]{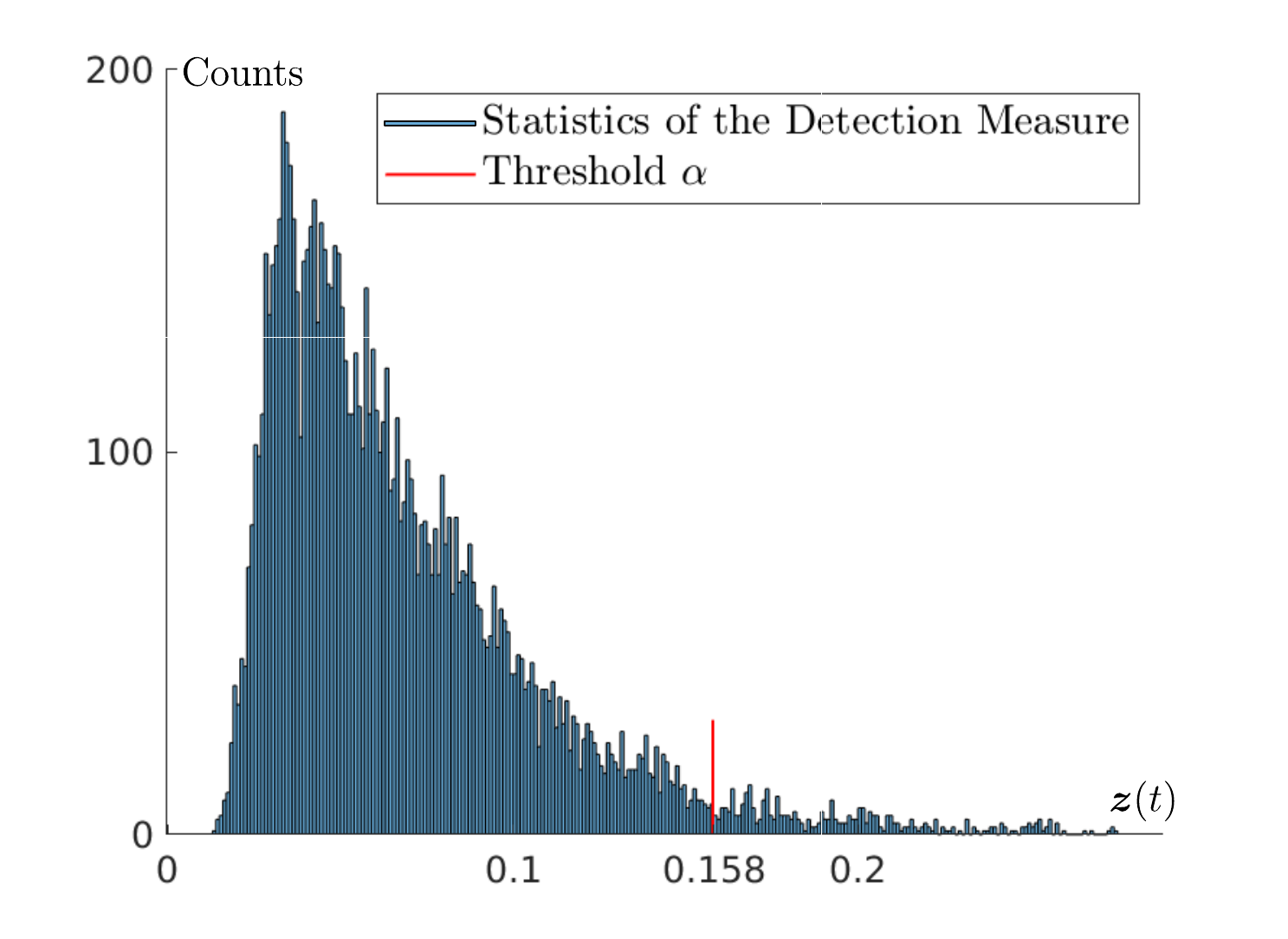}
         \caption{\small {Statistics of $z$.}}%
         \label{fig:z}%
     \end{minipage}
     \hfill{}
     \begin{minipage}[c]{0.64\columnwidth}
         \centering
         \includegraphics[width=\textwidth]{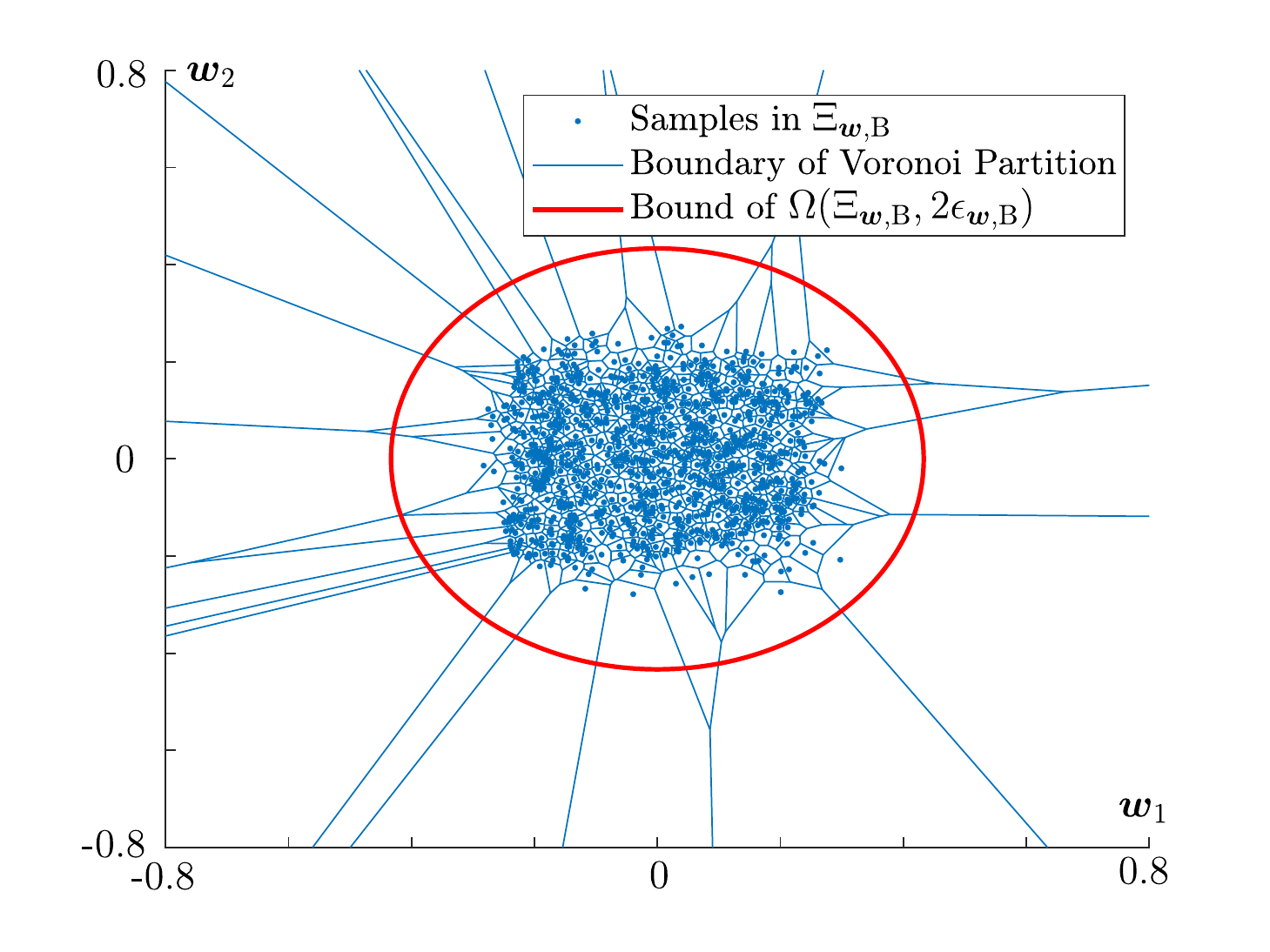}%
          \caption{\small {Probabilistic Support of $\prob_{\vect{w}}$.}}
         \label{fig:support}%
     \end{minipage}
     \hfill{}
     \begin{minipage}[r]{0.64\columnwidth}
         \centering
         \includegraphics[width=\textwidth]{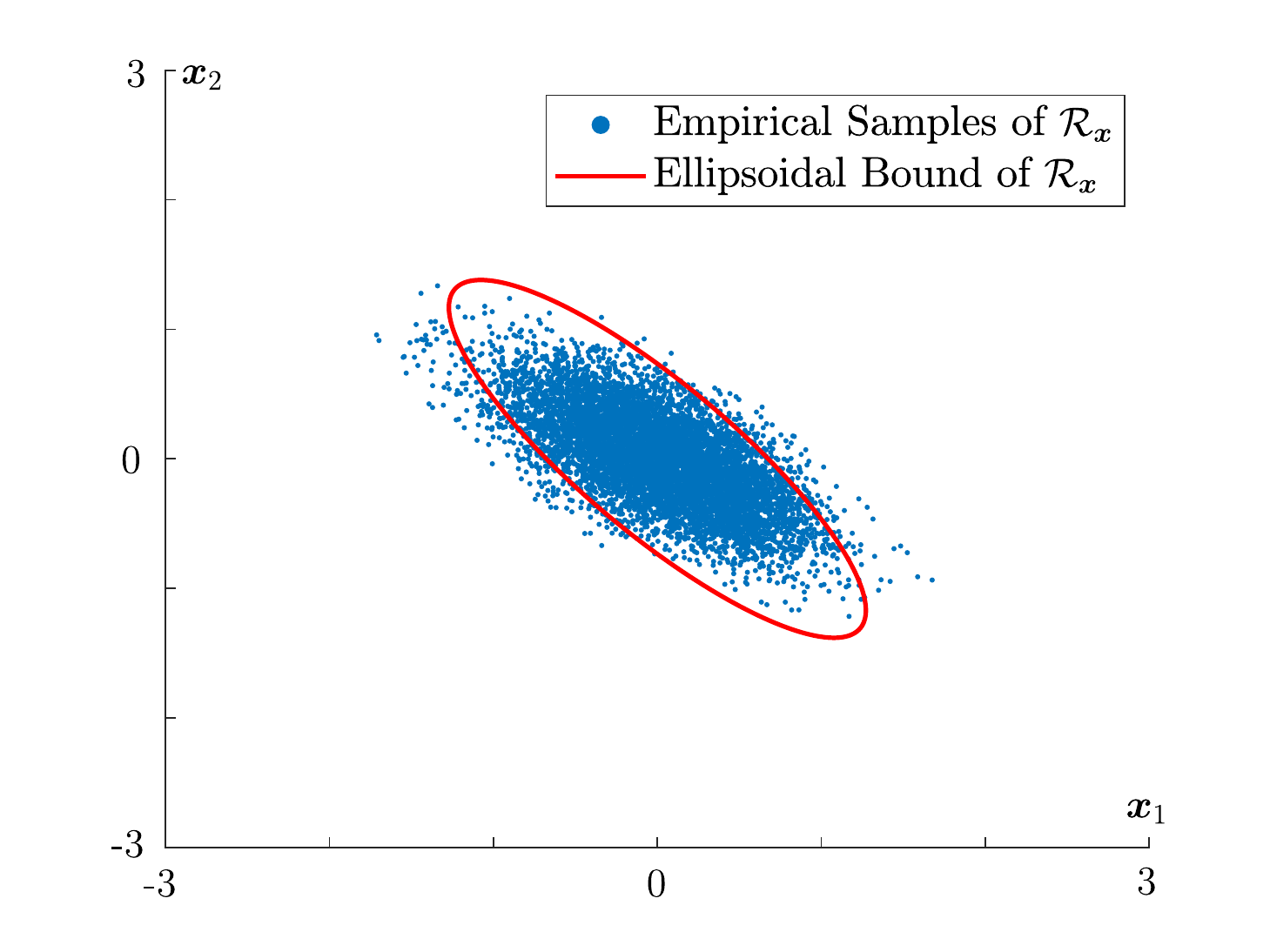}%
         \caption{\small {Empirical and Bound of $\mathcal{R}_{\vect{x}}$.}}%
         \label{fig:ReachableSet}%
     \end{minipage}
\end{figure*}
A tight reachable set bound can be now derived by solving
\begin{equation}
  \begin{aligned}
   \min\limits_{Q, a_1,a_2}  & \;  -\log \det(Q) , \\
       \st \; & \; \eqref{eq:Q}, \; \eqref{eq:W},
  \end{aligned} \label{eq:P2} \tag{P2}
\end{equation}
which is a convex semidefinite program, solvable via e.g.,
SeDuMi~\cite{JS:99}. Note that the probabilistic reachable set is
\begin{equation*}
  \mathcal{R}_{\vect{x}}:= \cup_{M=1}^{\infty} \mathcal{R}_{\vect{x},M},
\end{equation*}
which again can be approximated via $Q^{\star}$ solving~\eqref{eq:P2}
for\footnote{The set $\mathcal{R}_{\vect{x}}$ is in fact contained in
  the projection of $\mathcal{E}(Q^{\star})$ onto the state subspace,
  i.e., $\mathcal{R}_{\vect{x}} \subset \setdef{\vect{x}}{
    \trans{\vect{x}} \big( Q_{xx} - Q_{xe} Q_{ee}^{-1} \trans{Q_{xe}}
    \big) \vect{x} \leq \frac{(2-a)}{1-a}}$ with
  $Q^{\star}:=\begin{bmatrix} Q_{xx} & Q_{xe} \\ \trans{Q_{xe}} &
    Q_{ee}
\end{bmatrix}$. See, e.g.,~\cite{CM-IS-JR-DN:18} for details. }
\begin{equation*}
\mathcal{R}_{\vect{x}} \subset  \mathcal{E}(Q^{\star})=\setdef{\vect{x}\in\real^n}{ \trans{\vect{\xi}}Q^{\star}\vect{\xi} \leq \frac{(2-a)}{1-a}}.
\end{equation*}

\section{Simulations}\label{sec:simulation}
In this section, we demonstrate the performance of the proposed attack
detector, illustrating its distributional robustness w.r.t.~the
system noise. Then, we consider stealthy attacks as
in~\eqref{eq:gammabar} and analyze their impact by quantifying the
probabilistic reachable set and outer-approximation bound.

Consider the stochastic system~\eqref{eq:augment}, given as
\begin{equation*} \small
  \begin{aligned}
&    A= \begin{bmatrix}
     1.00 & 0.10 \\ -0.20 & 0.75
   \end{bmatrix} , \;
  B= \begin{bmatrix}
   0.10 \\ 0.20
  \end{bmatrix}, \;
  L=\begin{bmatrix}
    0.23  \\  -0.20
  \end{bmatrix}, \\
&  C=\begin{bmatrix}
   1 & 0
  \end{bmatrix}, \;
  K=\begin{bmatrix}
   -0.13 & 0.01
  \end{bmatrix}, \;
  n=2, \; m=p=1, \\
  & \vect{w}_1\sim \mathcal{N}(-0.25,0.02)+\mathcal{U}(0,0.5), \;  \vect{v} \sim \mathcal{U}(-0.3,0.3), \\
  & \vect{w}_2 \sim \mathcal{N}(0,0.04)+\mathcal{U}(-0.2,0.2), \; \\
  \end{aligned} \\
\end{equation*}
where $\mathcal{N}$ and $\mathcal{U}$ represent the normal and uniform
distributions, respectively. We consider $N=10^3$ benchmark samples
for $\prob_{\vect{r}, \Batch}$ and $T=10^2$ real-time samples for
$\prob_{\vect{r}, \Detect}$. We select the parameter $q=1$,
$\beta=0.01$ and false alarm rate $\Delta=0.05$. We select the prior
information of the system noise via parameters $a=1.5$, $c_1=1.84
\times 10^6$ and $c_2=12.5$. Using the measure-of-concentration
results, we determine the detector threshold to be $\alpha=0.158$. In
the normal system operation (no attack), we run the online detection
procedure for $10^4$ time steps and draw the distribution of the
computed detection measure $z(t)$ as in Fig.~\ref{fig:z}. We verify
that the false alarm rate is $3.68 \%$, within the required rate
$\Delta=5 \%$. When the system is subject to stealthy attacks, we
assume $\vect{\xi}_0=\vect{0}$ and visualize the Voronoi
partition $\mathcal{V}(\Xi_{\vect{w},\Batch})$ (convex sets with blue
boundaries) of the probabilistic support
$\Omega(\Xi_{\vect{w},\Batch},\epsilon_{\vect{w},\Batch})$ and its
estimated ellipsoidal bound (red line) as in
Fig.~\ref{fig:support}. Further, we demonstrate the impact of the
stealthy attacks~\eqref{eq:gammabar}, as in
Fig.~\ref{fig:ReachableSet}. We used $10^4$ empirical points of
$\mathcal{R}_{\vect{x}}$ as its estimate and provided an ellipsoidal
bound of $\mathcal{R}_{\vect{x}}$ computed by solution
of~\eqref{eq:P2}. It can be seen that the proposed probabilistic
reachable set effectively captures the reachable set in probability.
Due to the space limits, we omit the comparison of our approach to the
existing ones, such as the classical $\chi^2$ detector
in~\cite{YM-BS:15}
and the CUMSUM procedure~\cite{JM-DU-HS-KJ:18}.
However, the difference should be clear: our proposed approach is
robust w.r.t. noise distributions while others leverage the moment
information up to the second order, which only capture sufficient
information about certain noise distributions, e.g., Gaussian.

\section{Conclusions}\label{sec:Conclusion}
A novel detection measure was proposed to enable distributionally
robust detection of attacks w.r.t.~unknown, and light-tailed system
noise. The proposed detection measure restricts the behavior of the
stealthy attacks, whose impact was quantified via reachable-set
analysis.

 \bibliographystyle{IEEEtran}
\bibliography{../../bib/alias,../../bib/SMD-add,../../bib/JC,../../bib/SM}
\end{document}